\documentclass[11pt]{amsart}
\usepackage{graphicx}
\usepackage{amssymb,amsmath,amsfonts,amsthm,amscd,mathrsfs, epstopdf, amsbsy, psfrag, color}
\usepackage{wasysym}
\usepackage[all]{xy}
\usepackage{hyperref}
\input xy
\xyoption{all}
\CompileMatrices

\numberwithin{equation}{section}

\newtheorem{teo}[equation]{Theorem}
\newtheorem{lem}[equation]{Lemma}
\newtheorem{cor}[equation]{Corollary}
\newtheorem{prop}[equation]{Proposition}
\theoremstyle{definition}
\newtheorem{defi}[equation]{Definition}
\newtheorem{setup}[equation]{Setup}
\theoremstyle{remark}
\theoremstyle{Notation}
\newtheorem{nota}[equation]{Notation}

\theoremstyle{Delta's Hypothesis}
\newtheorem{hip}[equation]{Delta's Hypothesis}

\newenvironment{prova}{\noindent\textbf{Proof: }}{\hfill $\blacksquare$}
\newenvironment{provaprincipal}{\noindent\textbf{Proof [Theorem \ref{principal}]: }}{\hfill $\blacksquare$}

\newcommand{\R}{\mathbb{R}}

\newcommand{\Z}{\mathbb{Z}}
\newcommand{\N}{\mathbb{N}}

\newcommand{\p}{\mathbb{P}}

\newcommand{\fracc}{\operatorname{Frac}}

\newcommand{\ox}{\mathcal{O}_X}

\newcommand{\hp}{\operatorname{Hip}}
\newcommand{\kerr}{\operatorname{ker\,}}
\newcommand{\id}{\mathord{\rm id}}

\newtheorem{conjg}{Conjecture}

\newtheorem{conja}{Conjecture}

\newtheorem{conjga}{Conjecture}

\newif\ifproofmode
\proofmodetrue
\long\def\comment#1{\ifproofmode\bgroup
  \color{red} \tiny #1\egroup\fi}

\title[Free Structures in Division Rings]{Free Structures in Division Rings}

\author{R.~Fehlberg Junior}
\address{Universidade Federal do Esp\'irito Santo, Centro de Ci\^encias Exatas, Departamento de Matem\'atica, Vit\'oria, ES 29075-910,
  Brazil}
\email{fjrenato@yahoo.com.br}

\begin{document}

\maketitle

\begin{abstract}
Makar-Limanov's conjecture states that if a division ring $D$ is finitely generated and infinite dimensional over its center $k$ then $D$ contains a free $k$-subalgebra of rank $2$. In this work, we will investigate the existence of such structures in $D$, the division ring of fractions of the skew polynomial ring $L[t;\sigma]$, where $t$ is a variable and $\sigma$ is a $k$-automorphism of $L$. For instance, we prove Makar-Limanov's conjecture when either $L$ is the function field of an abelian variety or the function field of the $n$-dimensional projective space.
\end{abstract}

\section{Introduction}
Since the first example of a division ring was given (real quaternions: W.R. Hamilton in $1843$), this kind of object has fascinated a lot of mathematicians. Even though a classical structure, not so much is known about division rings, mainly those of infinite dimension over their center. One of the great contributions was given by A.I. Lichtman in $1977$, that conjectured in \cite{licht77} the following:

\begin{conjg}[Lichtman] The multiplicative group $D^\times$ of a division ring $D$ contains a noncyclic free group.
\end{conjg}

Some years later, L. Makar-Limanov proved the existence of free noncommutative algebras in the division ring of fractions of the first Weyl algebra over $\R$ \cite{ml83}. This breakthrough led Makar-Limanov to raise the following question in \cite{ml84algebra}:

\begin{conja}[Makar-Limanov]\label{conjml} Let $D$ be a division ring with center $k$. If $D$ is finitely generated and infinite dimensional over $k$ then $D$ contains a free $k$-subalgebra of rank $2$.
\end{conja}

Many authors have been studying the existence of free subalgebras \cite{ml83,ml84algebra,LOR,mlm91,GS,GS4,licht99,GT2,BR2}. Besides, some authors \cite{FGS, GT2, GS3, GS, JS} have been studying the following generalization of the aforementioned conjectures:

\begin{conjga}\label{conjga} Let $D$ be a division ring with center $k$. If $D$ is finitely generated and infinite dimensional over $k$ then $D$ contains a free group algebra over $k$.
\end{conjga}

Clearly conjecture GA implies in both conjectures A and G. In special cases \cite{licht84ii,GS3} conjecture A implies conjecture GA.

In the present text, we will focus our efforts on conjecture A.

Let $L/k$ be a field extension and let $D=L(t;\sigma)$ be the division ring of fractions of the skew polynomial ring $L[t;\sigma]$, where $\sigma$ is a $k$-automorphism of $L$. Our main result is the following:

\begin{teo}\label{principal} Let $L/k$ be a proper field extension such that $k$ is algebraically closed in $L$ and let $\sigma$ be a $k$-automorphism of $L$ such that $k=L^{\langle\sigma\rangle}$\footnote{The following consequence of a result by Makar-Limanov and Malcolmson (\cite{mlm91} Lemma 1) justifies our choices of $k$: a division ring $D$ over $k$ contains a free $k$-subalgebra of rank $2$ if and only if $D$ contains a free $K$-subalgebra of rank $2$ for any central subfield $K$.} $($fixed field by $\sigma$$)$. Consider $f\in A\setminus k$ for some $k$-subalgebra $A$ of $L$ such that Delta's Hypothesis \ref{hip} is true. Then, the elements $f$ and $t(1-t)^{-1}f$ generate a $k$-free subalgebra in $D=L(t;\sigma)$.
\end{teo}

\emph{Delta's Hypothesis} can be stated as follows:

\begin{hip}\label{hip} Let $L/k$ be a proper field extension such that $k$ is algebraically closed in $L$ and let $\sigma$ be a $k$-automorphism of $L$ such that $k=L^{\langle\sigma\rangle}$. Then

$$\Delta(L)\cap A\subset k,$$ where $\Delta=\sigma-\id$, for some $k$-subalgebra $A$ of $L$ such that $\fracc(A)=L$.

\end{hip}

This is a technical condition that will be only used in the proof of Lemma \ref{beta}, which is the main ingredient in the proof of Theorem \ref{principal}. For instance, if a proper field extension $L/k$ and a $k$-automorphism of $L$ satisfy Delta's Hypothesis \ref{hip}, then $L(t;\sigma)$ contains a free $k$-subalgebra of rank $2$. In Section 5, we will show cases in which this condition can be verified, that is, when either $L$ is the function field of an abelian variety or $L$ is the function field of the $n$-dimensional projective space over a field $k$, generalizing \cite{LOR,GT2}. Besides, applying a result by Lichtman \cite{licht84ii}, we will obtain a particular case of conjecture GA (see Corollary \ref{cor1}) from Theorem \ref{principal}, that emphasizes the importance of these techniques (see \cite{LOR,GT2}) and the generators found.

It would be remiss of us not to mention a recent work by Bell and Rogalski. In \cite{BR2}, these two authors proved, among other results, that if $L/k$ is a field extension with $k$ uncountable and $\sigma$ is a $k$-automorphism of $L$ then the existence of free subalgebras in $L(t;\sigma)$ is equivalent to $L$ having an element lying on an infinite $\sigma$-orbit. In the countable case, they assume that either $L/k$ is infinitely generated or $\sigma$ is induced by a $k$-automorphism of a quasi-projective variety. We emphasize our techniques are quite different from theirs and we do not impose restrictions on the base field, except in Section 5.

\section{Setup and definitions}
In this section we will define the main objects that will be used throughout the text.

\begin{setup}\label{setup1} Let
\begin{itemize}
\item $k$ be a field;
\item $k\subsetneq L$ be a field extension such that $k$ is algebraically closed in $L$;\footnote{Next, we will suppose that $L\otimes_kL$ is a domain. Then, we are assuming a little bit more about this field extension. For example, we will have this condition satisfied when $L/k$ is separable. For more details, consult \cite{SHC14},\cite{Lang} p.360-368 and \cite{stacksp} sections 7.41-7.47 and 28.4-28.7.}
\item $\sigma$ be a $k$-automorphism of $L$ such that $k=L^{\langle\sigma\rangle}$.
\end{itemize}
\end{setup}

\begin{nota} Assume Setup \ref{setup1}.
\begin{itemize}
\item We identify $L$ in $\fracc (L\otimes_k L)$ by $\ell\mapsto\ell\otimes 1$;
\item $\sigma_L=\id\otimes\sigma$ the automorphism of $L\otimes_k L$ induced by $\sigma$;
\item write again $\sigma_L$ as the automorphism of $\fracc (L\otimes_k L)$ induced by $\id\otimes\sigma$.
\end{itemize}
\end{nota}

\begin{defi} Assume Setup \ref{setup1}. We define:
\begin{enumerate}
\item $\mathscr{F}$ as the ring of functions $f:\mathbb{N}\rightarrow \fracc(L\otimes_k L)$ quotiented by the equivalence relation: $f\sim g$ iff $f(n)=g(n)$ for all $n\gg 0$. The sum and multiplication in $\mathscr{F}$ are coordinate to coordinate. Besides, $\fracc(L\otimes_k L)\cong\mathscr{C}\subset\mathscr{F}$, the set of constant functions in $\mathscr{F}$;
\item the shift operator $s:\mathscr{F}\rightarrow\mathscr{F}$ such that $\phi^s(n)=\phi(n+1)$, that is an automorphism of $\mathscr{F}$ such that restrict to $\mathscr{C}$ is the identity;
\item the operator $\Delta:\mathscr{F}\rightarrow\mathscr{F}$, given by $\Delta(\phi)=\phi^s-\phi$; using the same letter we define $\Delta$ analogously over $\fracc(L\otimes_k L)$ and $L$.
\end{enumerate}
\end{defi}

\section{Lemmas}

In this section we will show some lemmas, starting with results about the auxiliary ring $\mathscr{F}$. The hypotheses will be basically those in Setup \ref{setup1}. The proof of Lemma \ref{injet} is straightforward and it will be omitted.

\begin{lem}\label{injet} Assume Setup \ref{setup1}.
\begin{enumerate}
\item The $($left-$)$$L$-linear application
\begin{equation*}
    \begin{array}{cccc}
        \iota: & \fracc(L\otimes_kL) & \to & \mathscr{F} \\
           & \xi & \mapsto & (\sigma_L^j(\xi))_{j\geq 0} \\
    \end{array}
\end{equation*}
is a ring morphism $($therefore injective$)$.
\item The shift operator restricts to $\iota(\fracc(L\otimes_kL))$ and the following diagram is commutative.
$$\xymatrix{
& L\otimes_kL\ar@{^{(}->}[r]\ar@{^{(}->}[d]_{\sigma_L} & \fracc(L\otimes_kL) \ar@{^{(}->}[r]_{\hskip10mm\iota}\ar@{^{(}->}[d]_{\sigma_L} & \mathscr{F} \ar[d]^{s}\\
& L\otimes_kL\ar@{^{(}->}[r] & \fracc(L\otimes_kL) \ar@{^{(}->}[r]_{\hskip10mm\iota} & \mathscr{F}\\
}$$
\end{enumerate}
\end{lem}

The next lemma will show what happens with transcendental elements of $L$ over $k$ after a base change.

\begin{lem} Assume Setup \ref{setup1}. Consider $f\in L\setminus k$ $($therefore transcendental over $k$$)$. Then $1\otimes f$ is transcendental over $L$.
\end{lem}

\begin{prova}\label{trans2} If $L[1\otimes f]=L\otimes_k(k[f])$ then $\dim_LL[1\otimes f]=\dim_LL\otimes_k(k[f])=\dim_kk[f]$ and as $f$ is transcendental over $k$, we would have $1\otimes f$ transcendental over $L$. The proof of the first equality is simple.

\end{prova}

\begin{lem}\label{fixo} Assume Setup \ref{setup1}. Then $L^{\langle\sigma^{t}\rangle}$=k, where $L^{\langle\sigma^{t}\rangle}$ is the fixed field by $\sigma^t$ and $t\in\Z\setminus \{0\}$.
\end{lem}

\begin{prova}
First, if $\zeta\in L^{\langle\sigma^{t}\rangle}$ then $$\prod_{0\leq i<|t|}{(X-\sigma^i(\zeta))}$$ is a monic polynomial of degree $|t|$ canceled by $\zeta$ with coefficients in $k$. Therefore, the extension $L^{\langle\sigma^{t}\rangle}/k$ is algebraic. Since $k$ is algebraically closed in $L$, it follows that $L^{\langle\sigma^{t}\rangle}=k$.

\end{prova}

\begin{lem}\label{trans} Assume Setup \ref{setup1} and let $f\in L\setminus k$. Then, the constant function $\epsilon_f\stackrel{def}{=}1\otimes f\in\mathscr{C}\subset\mathscr{F}$ is transcendental over $\iota(\fracc(L\otimes_k L))$.
\end{lem}

\begin{prova}
Suppose that $$\sum_j{\xi_j(1\otimes f)^j}=0,$$ where $\xi_j\in\iota(\fracc(L\otimes_k L))$ aren't all nulls, that is, $\xi_j=(\sigma_L^n(e_j))_{n\geq 0}$ and $e_j\in\fracc(L\otimes_k L)$ for each $j$. Therefore,
\begin{equation}\label{eqlem}
\sum_j{[\sigma_L^n(e_j)](1\otimes f)^j}=0
\end{equation}
for all $n>n_0$, for some $n\gg 0$. With $n>n_0$ fixed, (\ref{eqlem}) is an equality in $\fracc(L\otimes_kL)$. Applying $\sigma_L^{-n}$ in (\ref{eqlem}) we obtain $$\sum_j{e_j[\sigma_L^{-n}(1\otimes f)^j]}=\sum_j{e_j[\sigma_L^{-n}(1\otimes f)]^j}=0$$ for each $n>n_0$. Therewith, if we consider $$g(X)=\sum_j{e_jX^j}\in(\fracc(L\otimes_k L))[X],$$ it follows that $g$ (a nonzero polynomial) would have $\sigma_L^{-n}(1\otimes f)$ as a root for all $n>n_0$. By Lemma \ref{fixo}, $1\otimes f$ isn't fixed for any power of  $\sigma_L$. In fact, if there is $r>0$ such that $\sigma^r_L(1\otimes f)=1\otimes f$, then $1\otimes\sigma^r(f)=1\otimes f$ we obtain $\sigma^r(f)-f=0$, that is a contradiction because Lemma \ref{fixo} and the choice of $f$. Therefore, the elements $\sigma_L^{-n}(1\otimes f)$ are all distinct for $n>n_0$ and $g(X)$ has infinite roots, that is an absurd.

\end{prova}

\begin{lem}\label{seqexata} Assume Setup \ref{setup1}. Let $A$ be a $k$-algebra. Then, $(A\otimes_kL)^{\langle\sigma_A\rangle}=A$, where $\sigma_A\stackrel{def}{=}\id\otimes\,\sigma$.
\end{lem}

\begin{prova} It is enough to apply $A\otimes_k-$ in the sequence $0\rightarrow k\rightarrow L\stackrel{\sigma-\id}{\rightarrow} L$.

\end{prova}

\section{Main result}\label{sectionmain}
Conditioned to Delta's Hypothesis, we will prove Theorem \ref{principal}. But before, we need to check Lemma \ref{beta} that is the main ingredient of its proof.

Because the generality degree of the statements and some specific modifications, we are including here the proofs of these results, but they are essentially the same as found in \cite{GT2}.

We start with some definitions:

\begin{defi}\label{defibeta} Assume Setup \ref{setup1}. Let $f\in L\setminus k$. We define:
\begin{enumerate}
\item $\gamma=\gamma_f\in\mathscr{F}$ by $\gamma(n)=(\id\otimes\sigma^n)(1\otimes f)=\sigma_L^n(1\otimes f)$ for all $n\in\N$;
\item a morphism of vector $k$-spaces (not rings!) $q:L[[t;\sigma]]\rightarrow\mathscr{F}$ by $\sum{t^na_n}\mapsto(1\otimes a_n)_{n\geq 0}$; \footnote{Note that $q(f)=(1\otimes f,0,0,\ldots)\notin\iota(\fracc(L\otimes L))$ and $\kerr(q)$ is formed by the polynomials in $L[[t;\sigma]]$.}
\item For $(n\in\mathbb{N})$ and $J=(j_1,j_2,\ldots,j_{s+1})\in\mathbb{N}^{s+1}_{>0}$ $($$s>0$$)$ put
$$\beta_J(n)=\sum_{n>n_1>\ldots>n_s>0}{\gamma(n_1)^{j_1}\gamma(n_2)^{j_2}\ldots\gamma(n_s)^{j_s}(1\otimes f)^{j_{s+1}}}.$$
For $s<0$ put $J=\emptyset$, $\beta_{\emptyset}(n)=1$ $\forall n>0$ and $\beta_\emptyset(0)=0$. For $s=0$, $\beta_J(n)=(1\otimes f)^{j_1}$ for all $n$, that is, $\beta_J=\epsilon^{j_1}\in\mathscr{C}$ is constant in $\mathscr{F}$.
\end{enumerate}
\end{defi}

The next result will be used only in the proof of Lemma \ref{beta}. Its proof is straightforward and it will be omitted.

\begin{lem}\label{lema3} Let $\psi,\varphi:\N\rightarrow\fracc(L\otimes_kL)$. Then:
\begin{enumerate}
\item $\Delta(\varphi\psi)=\varphi^{s}\Delta(\psi)+\Delta(\varphi)\psi$;
\item $\Delta(\beta_J)=\gamma^{j_1}\beta_{J(1)}$, where $J=(j_1,...,j_{s+1})\in\mathbb{N}^{s+1}_{>0}$, $s>0$ and $J(1):=(j_2,\ldots,j_{s+1})$. For $s=0$ we have $\Delta(\beta_J)=0$.
\end{enumerate}
\end{lem}

\begin{lem}\label{beta} Assume Setup \ref{setup1} and let $f\in A\setminus k$, for some $k$-subalgebra $A$ of $L$ that verifies Delta's Hypothesis \ref{hip}. Consider the elements $\beta_J$ as in Definition \ref{defibeta}. Then, the $\beta_J$ are linearly independent over $\iota(L_D)$ in $\mathscr{F}$.
\end{lem}

\begin{prova} Suppose the statement is not true. For a multi-index $J=(j_1,\cdots, j_s)$ write $|J|=s$ its size and if $\mathscr{J}$ is a finite set of index $J$, put $|\mathscr{J}|=\mbox{max}\{|J|; J\in\mathscr{J}\}$. We will prove by induction over $|\mathscr{J}|$.

Let
\begin{eqnarray}\label{soma}
\sum_{J\in\mathscr{J}}{e_J\beta_J}=0
\end{eqnarray}
be a nontrivial relation, where $e_J\in\iota(L_D)^\ast$ (always identifying with $\iota(e_J)$, for $e_J\in (L_D)^\ast$ using the same notation), chosen in the following way:

\begin{itemize}
\item[(i)] between all nontrivial relations over $\iota(L_D)$, choose that with $u=|\mathscr{J}|\geq 0$ minimal;
\item[(ii)] between the relations satisfying (i), choose that the number of $J$'s with $|J|=u$ is minimal.
\end{itemize}

If $u\leq 1$, we have $\beta_{(j_1)}=(1\otimes f)^{j_1}$ for $u=1$ e $\beta_\emptyset=1\in\mathscr{C}$ for $u=0$. Note that the $j_1$'s are distinct from each other. By Lemma \ref{trans}, these elements are linearly independent over $\iota(\fracc(L\otimes_kL))$ and therefore over $\iota(L_D)$.

For the general case ($u\geq 2$), rewrite (\ref{soma}) in the following form:

\begin{eqnarray}\label{eq1}
\sum_{|J|=u}{e_J\beta_J}=\sum_{|J|<u}{d_J\beta_J}
\end{eqnarray}

We can assume that $e_J=1$ for some $J$ such that $|J|=u$. Since $J=(j_1,\ldots,j_s)\in\N^s$, it follows from Lemma \ref{lema3} that
$$\Delta(e_J\beta_J)=\Delta (e_J)\cdot\beta_J+e_J^s\gamma^{j_1}\beta_{J(1)},$$

\begin{flushleft}
for $s\geq 2$ and $\Delta(e_J\beta_J)=\Delta (e_J)\cdot\beta_J$ if $s\leq 1$, where $\Delta(e_J), e_J^s\gamma^{j_1}\in\iota(L_D)$ ($L_D$ is stable by $\sigma_L$) and $|J(1)|<|J|$. In other words, applying $\Delta$ does not increase the number of $J$'s with maximum size.\end{flushleft}

Applying $\Delta$ in both sides of (\ref{eq1}) we obtain

\begin{eqnarray}\label{eq2}
\sum_{|J|=u}{\Delta (e_J)\cdot\beta_J}+\sum_{|J|=u}{e_J^s\gamma^{j_1}\beta_{J(1)}}=\Delta\left(\sum_{|J|<u}{d_J\beta_J}\right)
\end{eqnarray}

For $e_J=1$, we have $\Delta(e_J)=0$. Since $|J(1)|<|J|$, the left side of (\ref{eq2}) has at least one term less with $|J|=u$, while in the right side of (\ref{eq2}) all the terms have maximum size $|J|<u$. Since the coefficients in (\ref{eq2}) still are in $\iota(L_D)$, it follows from the minimality $(ii)$ that all of them are zero. Therefore, by the equality of terms with the same size of $J$, we obtain:

\begin{itemize}
  \item[(a)] $\Delta(e_J)=0$, for all $J$ such that $|J|=u$;
  \item[(b)] $\Delta(d_J)=\displaystyle\sum_{\stackrel{|I|=u}{I(1)=J}}{e_I^s\gamma^{i_1}}$ for all $J$ fixed with $|J|=u-1$, where $I=(i_1,\ldots,i_u)$.
\end{itemize}

Note that $$0=\Delta(e_J)=\iota(1\otimes\Delta(e_J))$$
and thus, $1\otimes\Delta(e_J)=0$ by the injectivity of $\iota$ and therefore $\Delta(e_J)=0$ ($e_J\in L$). Then, since $L^{\langle\sigma\rangle}=k$, $(a)$ occurs if, and only if, $e_J\in k^\ast$, for all $J$ with $|J|=u$. For $(b)$ we have:
$$\iota(\Delta(d_J))=\displaystyle\sum_{\stackrel{|I|=u}{I(1)=J}}{e_I^s\gamma^{i_1}}\stackrel{e_I\in k}{=} \displaystyle\sum_{\stackrel{|I|=u}{I(1)=J}}{e_I\gamma^{i_1}}=\iota(\displaystyle\sum_{\stackrel{|I|=u}{I(1)=J}}{e_I\otimes f^{i_1}})$$
and thus, by the injectivity of $\iota$ we have,

\begin{eqnarray}\label{eq3}
\Delta(d_J)=\displaystyle\sum_{\stackrel{|I|=u}{I(1)=J}}{e_I\otimes f^{i_1}}=\displaystyle\sum_{\stackrel{|I|=u}{I(1)=J}}{1\otimes(e_If^{i_1})}.
\end{eqnarray}

By (\ref{eq3}) it follows that

\begin{eqnarray}\label{eq4}
\Delta(d_J)=\displaystyle\sum_{\stackrel{|I|=u}{I(1)=J}}{e_If^{i_1}},
\end{eqnarray}
with $d_J\in L$. Then, $\Delta(d_J)\in\Delta(L)\cap A\subset k$ by \emph{Delta's Hypothesis} \ref{hip}.

Since for all $I$, each $i_1\geq 1$ and they are distinct from each other, it follows from the transcendence of $f$ over $k$ that the $e_I$'s in (\ref{eq4}) are null. Since this holds for each $J$ with $|J|=u-1$, we obtain $e_I=0$ for all $I$ with $|I|=u$, contradicting $(i)$.

\end{prova}

\begin{provaprincipal}
We will work initially in $L((t;\sigma))\supset D$. After, using the morphism $q$ and $\iota$, we will make the computations in $\mathscr{F}$. We will show that distinct monomials in $f$ and $t(1-t)^{-1}f$ are linearly independent over $k$. These monomials have the form

\begin{eqnarray*}
m_I & = & m_{(i_0,i_1,\cdots,i_v)}=f^{i_0}t(1-t)^{-1}f^{i_1}\cdots f^{i_{v-1}}t(1-t)^{-1}f^{i_v}\in L[[t;\sigma]]\\
 & = & \displaystyle{\sum_{n\geq 0}{t^n}}{f^{i_0\sigma^n}\displaystyle{\sum_{n>n_1>\cdots>n_{v-1}> 0}f^{i_1\sigma^{n_1}+\cdots+i_{v-1}\sigma^{n_{v-1}}+i_v}}}
\end{eqnarray*}
for $v>0$ and $m_I=f^{i_0}$ for $v=0$, where $v$ is the number of terms $t(1-t)^{-1}f$, $i_0\geq 0$ and $i_1,\cdots,i_v\geq 1$ (for more details, \cite{GT2} Proposition $9$).

Applying $q$ in each $m_I$ we will obtain $\overline{m_I}\stackrel{def}{=}q(m_I)=\gamma^{i_0}\beta_{I(1)}$ for $v>0$ and $q(m_{(i_0)})=(1\otimes f^{i_0},0,0,\ldots)$ for $v=0$.

Now, consider
\begin{eqnarray}\label{sum0}
\displaystyle{\sum_{I}{a_Im_I}}=\displaystyle{\sum_{I\in\mathscr{I}_0}{a_Im_I}}+\displaystyle{\sum_{I\notin\mathscr{I}_0}{a_Im_I}}=0,
\end{eqnarray}
where $a_I\in k$ and $\mathscr{I}_0=\{(i_0,i_1,\cdots,i_v);v=0\}$. Applying $q$ in (\ref{sum0}) we obtain:

$$0=\displaystyle{\sum_{I}{a_I\overline{m_I}}}=\displaystyle{\sum_{I\notin\mathscr{I}_0}{a_I\overline{m_I}}}=\displaystyle{\sum_{I\notin\mathscr{I}_0}{(a_I\gamma^{i_0})\beta_{I(1)}}}=\displaystyle{\sum_{J\neq\emptyset}}{\left(\sum_{I(1)=J}{a_I\gamma^{i_0}}\right)\beta_{J}},$$
where $a_I\gamma^{i_0}\in\iota(L_D)$.

By Lemma \ref{beta}, the $\beta_J$'s are linearly independent over $\iota(L_D)$, therefore, for each $J\neq\emptyset$, we have

$$0=\displaystyle{\sum_{I(1)=J}}{a_I\gamma^{i_0}}\stackrel{a_I\in k}{=}\iota\left(\displaystyle{\sum_{I(1)=J}}{a_I\otimes f^{i_0}}\right).$$

Thus, $$0=\displaystyle{\sum_{I(1)=J}}{a_I\otimes f^{i_0}}\stackrel{a_I\in k}{=}1\otimes\left(\displaystyle{\sum_{I(1)=J}}{a_If^{i_0}}\right)$$ by the injectivity of $\iota$.

Therefore $$\displaystyle{\sum_{I(1)=J}}{a_If^{i_0}}=0.$$

Since the $i_0$'s are distinct one of each other, it follows from the transcendence of $f$ over $k$ that $a_I=0$ for all $I\notin\mathscr{I}_0$. Expression (\ref{sum0}) reduces to
$$\displaystyle{\sum_{I\in\mathscr{I}_0}{a_Im_I}}\stackrel{def}{=}\displaystyle{\sum_{I=(i_0)\in\mathscr{I}_0}{a_If^{i_0}}}=0,$$
and we apply again the transcendence of $f$ over $k$ to obtain $a_I=0$ for all $I\in\mathscr{I}_0$. This ends the proof.

\end{provaprincipal}

As a consequence of Theorem \ref{principal} we will obtain a particular case of conjecture GA.

\begin{cor}\label{cor1} Assume Setup \ref{setup1}. Let $f\in A\setminus k$ for some $k$-subalgebra $A$ of $L$ such that Delta's Hypothesis \ref{hip} is true. Then, $D=L(t;\sigma)$ contains a free group algebra over $k$.
\end{cor}

\begin{prova} By Theorem \ref{principal}, the elements $f$ and $g=t(1-t)^{-1}f$ generate a free $k$-subalgebra in $D=L(t;\sigma)$. Therefore the elements $fg$ and $g$ still generate a free $k$-subalgebra in $D$.

Consider in $L((t;\sigma))$ a valuation $v:L((t;\sigma))\mapsto \Z\cup\{\infty\}$ such that: $v$ is trivial over $L$ and $v(t)=1$. Thus, $$v(fg)=v(g)=1.$$

By \cite{licht84ii} p.524 Corollary 1, the elements $1+fg$ and $1+g$ generate a free group algebra in $D$ over $k$.

\end{prova}

\section{Applications}
In this section we will prove \emph{Delta's Hypothesis} in two cases: when either $L$ is the function field of an abelian variety or the function field of the $n$-dimensional projective space over a field $k$. We start with a lemma.

\begin{lem}\label{lbase} Let $k$ be a field and let $X$ be a regular separated integral noetherian scheme over $k$, with function field $L$. Let $\Sigma$ be an automorphism of $X$ and $\sigma$ the automorphism of $L$ induced by $\Sigma$. Then, if $\Delta(g)$ has at most one pole (not counting multiplicities), it follows that $g$ and $\sigma(g)$ have the same set of poles without counting multiplicities.
\end{lem}

\begin{prova} Let $g\in L$ be a nonzero element and let $S\subset X$ be the finite set of poles of $g$ without counting multiplicities.

We have $T=\Sigma^\ast S=\{\Sigma^{-1}D; D\in S\}$ is the set of poles of $\sigma(g)$. Therefore, $|S|=|T|$. The set of poles of $\Delta(g)$ contains the poles of $g$ and $\sigma(g)$, possible excluding the elements of $S\cap T$, that is, it contains $S\cup T\setminus S\cap T$. This set has cardinality $2(|S|-|S\cap T|)$.

Suppose $\Delta(g)$ has at most one pole (not counting multiplicities). Since $\Delta(g)\supset S\cup T\setminus S\cap T$ we obtain $|S\cup T\setminus S\cap T|\leq 1$. But $|S\cup T\setminus S\cap T|=2(|S|-|S\cap T|)$ and therefore $|S\cup T\setminus S\cap T|=0$, that is, $S=T$.

\end{prova}

The next lemma allows us to assume that the field $k$ is algebraically closed.

\begin{lem} Suppose that $\overline{\Delta}(\overline{k}\otimes_kL)\cap \overline{k}\otimes_kA\subset \overline{k}$, where $\overline{\Delta}=\id\otimes\Delta$. Then, $\Delta(L)\cap A\subset k$.
\end{lem}

\begin{prova} Let $f\in L$ and suppose $\Delta(f)\in A$. We want to show that $\Delta(f)\in k$. Since $\Delta(f)\in A$, it follows that $1\otimes\Delta(f)\in\overline{k}\otimes_kA$. But, $1\otimes\Delta(f)=\overline{\Delta}(1\otimes f)$. Then, by hypothesis, $1\otimes\Delta(f)\in\overline{k}=(\overline{k}\otimes_kL)^{\langle\id\otimes\sigma\rangle}$, by Lemma \ref{seqexata}. Therefore, $\Delta(f)\in k=L^{\langle\sigma\rangle}$.

\end{prova}

In order to prove Delta's Hypothesis, we will provide a necessary condition, at least when $X$ is a projective variety.

\begin{prop}[\emph{Delta's Hypothesis}]  Let $(X,\ox)$ be a projective variety over an algebraically closed field $k$ with function field $L$, $\Sigma$ an automorphism of $X$ and $D$ a prime divisor. Let $\sigma$ be the automorphism of $L$ induced by $\Sigma$. Assume the following condition:

\begin{equation}\label{eqd}
(\Sigma^n)^\ast(D)\neq D,
\end{equation}
for all integer $n\geq 1$. Then
$$\Delta(L)\cap\ox(X\setminus D)\subset k.$$
\end{prop}

\begin{prova} Let $g\in L$ and let $S$ be the set of poles of $g$ not counting multiplicities. Let also $T$ be the set of poles of $\sigma(g)$.
Suppose that $\Delta(g)\in\ox(X\setminus D)$. Any element of $\ox(X\setminus D)$ has at most one pole in $D$. By Lemma \ref{lbase}, it follows that $S=T=\Sigma^\ast S$. However, if $S\neq\emptyset$ then the set of poles of $g$ would be stabilized by $\Sigma$ and therefore, some power of $\Sigma$ would fix all the elements of $S$. By hypothesis, $D\notin S=T$ and then $D$ is not a pole of $\Delta(g)$, that must be constant.

\end{prova}

Finally, we will show that we can assume equation (\ref{eqd}) when either $X$ is an abelian variety or the $n$-dimensional projective space. The cases that we will consider are described by the following two setups.

\begin{setup}\label{setup3} Let
\begin{itemize}
\item $k$ be an algebraically closed field (which is not the algebraic closure of a finite field);
\item $X=(X,\ox)$ be an abelian variety over $k$ with $g=\dim_k X\geq 2$ $($see \cite{GT2} for the $1$-dimensional case$)$;
\item $L$ be the function field of $X$;
\item $H$ be a generic hyperplane (Bertini's Theorem: \cite{HAR} p.179).
\end{itemize}
\end{setup}

\begin{setup}\label{setup4} Let
\begin{itemize}
\item $k$ be a field;
\item $X=\p^n_k=(X,\ox)$
\item $L$ the function field of $X$;
\item $\Sigma$ be an automorphism of $X$ of infinite order such that the induced automorphism $\sigma$ satisfies $k=L^{\langle\sigma\rangle}$.
\end{itemize}
\end{setup}

\begin{lem}\label{lematrans} Assume Setup \ref{setup3}. There exists a point $P_0\in X$ such that the translation by $P_0$, $\Sigma_{P_0}=\Sigma:X\rightarrow X$, satisfies $(\Sigma^n)^\ast(D)\neq D$ for all integer $n\geq 1$, where $D=X\cap H$.

\end{lem}

\begin{prova}
Since abelian varieties are projective (\cite{milne} Theorem 6.4 p.29),
it follows that $H_n=\{x\in X; (\Sigma_x^n)^\ast(D)= D\}< X$ is finite for all integer $n\geq 1$ (\cite{MUN} Application 1 p.57),
where $\Sigma_x$ is the translation by $x$. Besides, $H_n=\{x\in X; nx\in H_1\}$ since $(\Sigma^n_x)^\ast=\Sigma^\ast_{nx}$. Then, $|H_n|=n^{2g}|H_1|$ where $n^{2g}=|K_n|=|\{x\in X; nx=0\}|$. Because of this, we obtain all the elements of $H_n$ from $H_1$ for all $n\geq 2$. Therefore, to obtain $x\in X$ such that $(\Sigma_x^n)^\ast(D)\neq D$ for all $n\geq 1$, it is enough to find $x\in X$, such that $x$ is $\Z$-linearly independent of the elements of $H_1$ (that will have infinite orbit). But, the rank of $X(k)$ (set of $k$-rational points) is infinite (\cite{FJ} Theorem 10.1 p.126), which concludes the proof.

\end{prova}

\begin{lem}\label{lemaproj1} Assume Setup \ref{setup4}. There exists a hyperplane $H$ of $\p^n_k$ such that $\Sigma^iH\neq H$ for all $i>0$.
\end{lem}

\begin{prova}
We have a bijection

$$\begin{array}{cccc}
\Gamma: &\hp(\p^n_k) & \rightarrow & (\p^n_k)^\ast\\
 & H=\{h=0\} & \mapsto & [h]
\end{array}$$
between the set of hyperplanes of $\p^n_k$ and the points of $(\p^n_k)^\ast$, where $h\in (k^{n+1})^\ast$ and $[h]$ is the class of $h$ in $(\p^n_k)^\ast$. Besides, $\Sigma$ induces an automorphism

$$\begin{array}{cccc}
\Sigma': &\hp(\p^n_k) & \rightarrow & \hp(\p^n_k)\\
 & H=\{h=0\} & \mapsto & \{\Sigma^\sharp(h)=h\circ\Sigma^{-1}=0\}
\end{array}$$

Note that $h\circ\Sigma^{-1}(p)=0 \Leftrightarrow\Sigma^{-1}(p)\in H\Leftrightarrow p\in\Sigma(H)$. Thereby, we have an automorphism $\Sigma'':(\p^n_k)^\ast\rightarrow(\p^n_k)^\ast$, defined by $[h]\mapsto [h\circ\Sigma^{-1}]$. Therefore, we have the following commutative diagram

$$\xymatrix{
& \hp(\p^n_k) \ar@{^{(}->}[r]_{\Gamma}\ar@{^{(}->}[d]_{\Sigma'} & (\p^n_k)^\ast \ar[d]^{\Sigma''}\\
& \hp(\p^n_k) \ar@{^{(}->}[r]_{\Gamma} & (\p^n_k)^\ast\\
}$$

This gives us a bijective correspondence between automorphisms of $\hp(\p^n_k)$ and automorphisms of $(\p^n_k)^\ast$.

Consider $M=\{H; H\,\, \mbox{is a hyperplane such that} \,\,\Sigma^iH=H\,\,\,\mbox{for some}\, i>0\}$. Let $H_i=\{h_i=0\}$, where $i=1,\ldots,n+2$, be distinct hyperplanes in $M$ fixed by powers $\alpha_i$ of $\Sigma$ respectively, with $i=1,\ldots,n+2$. Thereby $(\Sigma')^{\alpha_i}(H_i)=H_i$. By the diagram, these hyperplanes correspond to distinct points $P_i$ in $(\p^n_k)^\ast$, that are fixed by the respective powers of $\Sigma''$. Then, $(\Sigma'')^{\alpha_1\alpha_2\ldots\alpha_{n+2}}=(\Sigma'')^\alpha$ fix $n+2$ points and by a general fact (if an element of $\operatorname{PGL}_k(n)$ fix $n+2$ points of $\p^n_k$ then it is the identity) it follows that $(\Sigma'')^\alpha=\id$ implying that $(\Sigma')^\alpha=\id$. Therefore, $\Sigma^\alpha$ fix all hyperplane of $\p^n_k$. It is known that all point $x\in\p^n_k$ can be viewed as the intersection of $n$ hyperplanes, then $\Sigma^\alpha(x)=x$ for all $x\in\p^n_k$. Thus, $\Sigma^\alpha=\id$, that is an absurd. Therefore, $|M|\leq n+1$. By counting the elements of $(\p^n_k)^\ast$, there exists a point $P$ that it is not fixed by any power of $\Sigma''$, that corresponds to a hyperplane $H$ of $\p^n_k$ that it is not fixed by any power of $\Sigma$.

\end{prova}

\section*{Acknowledgments}
The author would like to thank his PhD advisors Daniel Levcovitz and Eduardo Tengan for all the help and outstanding support. The author was supported by FAPESP Grant 2009/54547-5 during his PhD.

\bibliographystyle{alpha}
\bibliography{fjrenato}

\def\cprime{$'$}
\begin{thebibliography}{MLM91}

\bibitem[BR12]{BR2}
J.~P. Bell and D.~Rogalski.
\newblock Free subalgebras of quotient rings of ore extensions.
\newblock {\em Algebra Number Theory}, 6(7):1349--1367, 2012.

\bibitem[Car]{SHC14}
P.~Cartier.
\newblock {\em S\'eminaire {H}enri {C}artan de l'{E}cole {N}ormale
  {S}up\'erieure, 1955/1956. {E}xtensions r\'eguli\`eres}.
\newblock tome 8, exp.14: 1-10.

\bibitem[FGS96]{FGS}
L.~M.~V. Figueiredo, J.~Z. Gon{\c{c}}alves, and M.~Shirvani.
\newblock Free group algebras in certain division rings.
\newblock {\em J. Algebra}, 185(2):298--313, 1996.

\bibitem[FJ74]{FJ}
G.~Frey and M.~Jarden.
\newblock Approximation theory and the rank of abelian varieties over large
  algebraic fields.
\newblock {\em Proc. London Math. Soc. (3)}, 28:112--128, 1974.

\bibitem[GT12]{GT2}
J.~Z. Gon{\c{c}}alves and E.~Tengan.
\newblock Free group algebras in division rings.
\newblock {\em Internat. J. Algebra Comput.}, 22(5), 2012.

\bibitem[Har77]{HAR}
R.~Hartshorne.
\newblock {\em Algebraic geometry}.
\newblock Springer-Verlag, New York, 1977.
\newblock Graduate Texts in Mathematics, No. 52.

\bibitem[Lan02]{Lang}
S.~Lang.
\newblock {\em Algebra}, volume 211 of {\em Graduate Texts in Mathematics}.
\newblock Springer-Verlag, New York, third edition, 2002.

\bibitem[Lic77]{licht77}
A.~I. Lichtman.
\newblock On subgroups of the multiplicative group of skew fields.
\newblock {\em Proc. Amer. Math. Soc.}, 63(1):15--16, 1977.

\bibitem[Lic84]{licht84ii}
A.~I. Lichtman.
\newblock On matrix rings and linear groups over fields of fractions of group
  rings and enveloping algebras. {II}.
\newblock {\em J. Algebra}, 90(2):516--527, 1984.

\bibitem[Lic99]{licht99}
A.~I. Lichtman.
\newblock Free subalgebras in division rings generated by universal enveloping
  algebras.
\newblock {\em Algebra Colloq.}, 6(2):145--153, 1999.

\bibitem[Lor86]{LOR}
M.~Lorenz.
\newblock On free subalgebras of certain division algebras.
\newblock {\em Proc. Amer. Math. Soc.}, 98(3):401--405, 1986.

\bibitem[Mil]{milne}
J.S. Milne.
\newblock Abelian varieties.
\newblock {\em http://www.jmilne.org/math/CourseNotes/AV.pdf}.

\bibitem[ML83]{ml83}
L.~Makar-Limanov.
\newblock The skew field of fractions of the {W}eyl algebra contains a free
  noncommutative subalgebra.
\newblock {\em Comm. Algebra}, 11(17):2003--2006, 1983.

\bibitem[ML84]{ml84algebra}
L.~Makar-Limanov.
\newblock On free subobjects of skew fields.
\newblock In {\em Methods in ring theory ({A}ntwerp, 1983)}, volume 129 of {\em
  NATO Adv. Sci. Inst. Ser. C Math. Phys. Sci.}, pages 281--285. Reidel,
  Dordrecht, 1984.

\bibitem[MLM91]{mlm91}
L.~Makar-Limanov and P.~Malcolmson.
\newblock Free subalgebras of enveloping fields.
\newblock {\em Proc. Amer. Math. Soc.}, 111(2):315--322, 1991.

\bibitem[Mum08]{MUN}
D.~Mumford.
\newblock {\em Abelian varieties}, volume~5 of {\em Tata Institute of
  Fundamental Research Studies in Mathematics}.
\newblock Published for the Tata Institute of Fundamental Research, Bombay,
  2008.
\newblock With appendices by C. P. Ramanujam and Yuri Manin, Corrected reprint
  of the second (1974) edition.

\bibitem[S{\'{a}}n]{JS}
J.~S{\'{a}}nchez.
\newblock Free group algebras in {M}alcev-{N}eumann skew field of fractions.
\newblock {\em {\tt arXiv:1107.2429 [math.RA]}}.

\bibitem[SG96]{GS3}
M.~Shirvani and J.~Z. Gon{\c{c}}alves.
\newblock On free group algebras in division rings with uncountable center.
\newblock {\em Proc. Amer. Math. Soc.}, 124(3):685--687, 1996.

\bibitem[SG98]{GS}
M.~Shirvani and J.~Z. Gon{\c{c}}alves.
\newblock Free group algebras in the field of fractions of differential
  polynomial rings and enveloping algebras.
\newblock {\em J. Algebra}, 204(2):372--385, 1998.

\bibitem[SG99]{GS4}
M.~Shirvani and J.~Z. Gon{\c{c}}alves.
\newblock Large free algebras in the ring of fractions of skew polynomial
  rings.
\newblock {\em J. of the London Math. Soc. (2)}, 60(2):481--489, 1999.

\bibitem[sta]{stacksp}
Stacks project.
\newblock {\em http://stacks.math.columbia.edu}.

\end{thebibliography}

\end{document}